\documentclass[12pt]{article}
\usepackage{amsmath}
\usepackage{amssymb}
\usepackage[cp1251]{inputenc}
\usepackage[russian]{babel}
\newcommand{\il}[2]{\int\limits_{#1}^{#2}}
\newcommand{\ilp}[1]{\int\limits_{#1}^{+\infty}}

\newcommand{\ph}{\phantom{a}}
\newcommand{\phh}{\phantom{aaa}}

\newcommand{\sist}[2]{\left\{
\begin{array}{l}
{#1}\\
\ph\\
{#2}
\end{array}
\right.}
\newcommand{\mb}[1]{\ph\mbox{#1}\ph}

\newcommand{\matr}[2]{\left (
\begin{array}{l}
{#1}\\
{#2}
\end{array}
\right)}

\newcommand{\matrr}[3]{\left (
\begin{array}{l}
{#1}\\
{#2}
\ph\\
{#3}
\end{array}
\right)}

\newcommand{\matrrr}[3]{\left (
\begin{array}{l}
{#1}\\
{#2}\\
{#3}
\end{array}
\right)}

\begin{document}

MSC 34C10

\vskip 20pt

\centerline{\bf  Oscillatory criteria for the hamiltonian systems}

\vskip 10 pt

\centerline{\bf G. A. Grigorian}

\vskip 10 pt

\centerline{0019 Armenia c. Yerevan, str. M. Bagramian 24/5}
\centerline{Institute of Mathematics NAS of Armenia}
\centerline{E - mail: mathphys2@instmath.sci.am, \ph phone: 098 62 03 05, \ph 010 35 48 61}

\vskip 20 pt

\noindent
Abstract. The Riccati equation method is used to establish some new oscillatory criteria for the hamiltonian systems in a new direction, which  is to break the positive definiteness restriction imposed on one of  coefficients of the hamiltonian system. The obtained results are compared with some known oscillatory criteria.
\vskip 20 pt


\vskip 20 pt

\centerline{\bf \S 1. Introduction}.

\vskip 10 pt

Let $A(t), \ph B(t)$ and $C(t)$ be complex valued continuous matrix - functions of dimension $n \times n$ on $[t_0;+\infty)$ and let $B(t)$ and $C(t)$ be hermitian matrices,  i.e. \linebreak $B(t) = B^*(t), \ph C(t) = C^*(t), \ph t\ge t_0$, where $*$ is the conjugation sign. Consider the hamiltonian system
$$
\sist{\Phi' = A(t) \Phi + B(t) \Psi;}{\Psi' = C(t)\Phi - A^*(t) \Psi, \phh t\ge t_0.} \eqno (1.1)
$$
Here $\Phi$ and $\Psi$ are the unknown  continuously differentiable on $[t_0;+\infty)$ matrix - functions of dimension  $n \times n$.

{\bf Definition 1.1}. {\it A solution $(\Phi(t), \Psi(t))$ of the system (1.1) is called prepared (or preferred), if $\Phi^*(t) \Psi(t) = \Psi^*(t) \Phi(t), \ph t\ge t_0.$}

{\bf Definition 1.2}. {\it A prepared solution $(\Phi(t), \Psi(t))$ of the system (1.1) is called oscillatory if $\det \Phi(t)$ has arbitrary large zeroes}.

{\bf Definition 1.3}. {\it The system (1.1) is called oscillatory if its every prepared solution is oscillatory.}

{\bf Definition 1.4}. {\it  A prepared solution $(\Phi(t), \Psi(t))$  of the system (1.1) is called oscillatory on the interval $[a;b] \ph (\subset [t_0;+\infty))$ if $\det \Phi(t)$ vanishes on $[a;b]$.}

{\bf Definition 1.5}. {\it The system (1.1) is called oscillatory on the interval $[a;b] \ph (\subset [t_0;+\infty))$ if its every prepared solution is oscillatory on $[a;b]$.}

Study of the  oscillatory behavior of the system (1.1) is an important problem of the qualitative theory of differential equations and many works are devoted to it (see [1 - 4] and cited works therein). In the works [1] and [2] some oscillatory criteria are proved for the system (1.1)  in terms of the coefficients $B(t)$ and $C(t)$ and  the  fundamental matrix  of the linear system $v' = A(t) v, \ph t\ge t_0$. In the works [3] and [4] some oscillatory criteria are obtained for the system (1.1) in terms of its coefficients. In all these criteria the positive definiteness condition on $[t_0;+\infty)$  is imposed on the coefficient $B(t)$ of the system (1.1) (therefore $B(t)$ is invertible for all $t\ge t_0$). The goal  of this paper is to obtain some oscillatory criteria for the system (1.1) in a new direction, which is to break the positive definiteness  restriction imposed on $B(t)$   for all $t\ge t_0$.

\vskip 20pt

\centerline{\bf \S 2.  Oscillatory criteria}

\vskip 20pt

{\bf 2.1. Main results}. Let $a_{jk}(t), \ph t\ge t_0, \ph j,k =1,2,$ be real valued continuous functions on $[t_0;+\infty)$. Along with the system (1.1) consider the following scalar one
$$
\sist{\phi' = a_{11}(t) \phi + a_{12}(t) \psi;}{\psi' = a_{21}(t) \phi + a_{22}(t) \psi, \ph t \ge t_0.} \eqno (2.1)
$$

{\bf Definition 2.1}. {\it The  system (2.1) is called oscillatory (on the interval $[a;b] \linebreak (\subset [t_0;+\infty))$) if for its every solution $(\phi(t), \psi(t))$ the function $\phi(t)$ has arbitrary large zeroes (the function $\phi(t)$ vanishes on $[a;b]$)}.

Hereafter we will assume that $[a;b]$ is a interval from the set $[t_0; +\infty)$.
Denote by $\Omega_n$ the set of all normal matrices $M$ (i. e. $M^* M = M M^*$) of dimension $n\times n$   each of which have eigenvalues $\lambda_1 = \lambda_1(M), \dots, \lambda_n = \lambda_n(M)$  with $Re \lambda_1 = \dots = Re \lambda _n \stackrel{def}{=}~ W(M)$. Note that  $M\in \Omega_n$ if in particular $M = \alpha I + i H$, where $\alpha \in (-\infty; +\infty), \ph I$  and $H$ are the identity and a hermitian matrices of dimension  $n\times n$ respectively. Denote by $\lambda(H)$ the least eigenvalue of any hermitian matrix $H$, and the non negative (positive) definiteness of  $H$ we denote by the symbol $H \ge 0 \ph (> 0)$. The trace of arbitrary square matrix M we denote by $tr (M)$. For any continuously differentiable hermitian matrix - function $S(t)$ set:\linebreak  $D_S(t) \equiv S'(t) + S(t) B(t) S(t) + A^*(t) S(t) + S(t) A(t) - C(t), \ph \sigma_S(t) \equiv W(A^*(t) + S(t) B(t))$.

\pagebreak

{\bf Theorem 2.1.} {\it Let the following conditions be satisfied:

\noindent
1) \ph $B(t) \ge 0, \ph t \in [a;b];$

\noindent
2) \ph there exists a continuously differentiable hermitian matrix - function $S(t)$  on $[a;b]$ such that $A^*(t) + S(t) B(t) \in \Omega _n, \ph t\in [a;b]$;

\noindent
3) the scalar system
$$
\sist{\phi' = \phh \sigma_S(t) \phi \ph +\ph  \frac{\lambda(B(t))}{n} \psi;}{\psi' = - tr(D_S(t)) \phi\ph - \ph \sigma_S(t) \psi, \ph t\in [a;b],} \eqno (2.2)
$$
is oscillatory on $[a;b]$.

\noindent
Then the system (1) is also oscillatory on $[a;b]$.} $\Box$

Indicate some particular cases in which the condition 2) of Theorem 2.1 is satisfied:

\noindent
$I) \ph A^*(t) \in \Omega_n, \ph t \ge t_0 \ph (S(t) \equiv 0);$

\noindent
$II) \ph A^*(t) = A_1(t) + \matr{A_2(t) \ph 0}{0 \ph A_3(t)}, \ph B(t) = \matr{B_1(t) \ph 0}{0\ph B_2(t)}, \ph t \in [a;b],$ where $A_2(t)$ and $B_1(t)$ are some matrix - functions of dimension $m\times m \ph (m< n), \ph A_1(t) \in \Omega_n,  \ph A_3(t) \in \Omega_{n -m}, \linebreak  W(A_3(t)) \equiv~ 0, \ph B_1(t) > 0, \ph t\in [a;b]:$

\noindent
$II_1) \ph A_2(t) B_1(t) = B_1(t) A_2(t), \ph t\in [a;b], \ph A_2(t) B_1^{-1}(t)$ is a continuously differentiable \linebreak matrix - function on $[a;b], \ph \biggl(S(t) \equiv \matr{- A_2(t) B_1^{-1}(t) \ph \hskip 3pt 0}{\phh 0 \phh \phh \phh 0}, \ph t \in [a;b]\biggr);$

\noindent
$II_2) \ph [A_2(t) + A_2^*(t)] B_1(t) = B_1(t)[A_2(t) + A_2^*(t)], \ph t\ge t_0, \ph [A_2(t) + A_2^*(t)] B_1^{-1}(t)$ is a continuously differentiable matrix - function on $[a;b];$
$$
\left( S(t) \equiv \matr{-\frac{A_2(t) + A_2^*(t)}{2} B_1^{-1}(t) \phh \hskip -1pt 0}{\phh 0 \phh \phh \phh \phh \phantom{aa}0}, \ph t \in [a;b] \right);
$$

\noindent
$III) \ph A^*(t) = A_1(t) + a(t) J, \ph B(t) = b(t) J, \ph J^2 = J = J^* = const \ph \Biggl(e.g. \ph J=\frac{1}{n}\matrrr{1...1}{.....}{1...1}\Biggr), \linebreak A_1(t) \in \Omega_n, \ph t \in [a;b],$ where $a(t)$ and $b(t)$ are some real valued continuous functions on $[a;b], \ph b(t) \ge ~0,\ph t\in [a;b], \ph \frac{a(t)}{b(t)}$  is a continuously differentiable function on $[a;b] \linebreak (S(t) \equiv - \frac{a(t)}{b(t)} J);$

\noindent
$IV) \ph A^*(t) = A_1(t) B(t), \ph t\in [a;b], \ph A_1(t) = A_1^*(t), \ph t\in [a;b],$ where $A_1(t)$ is a continuously differentiable matrix - function on $[a;b] \ph (S(t) \equiv - A_1(t), \ph t\in [a;b])$.

{\bf Remark 2.1.} {\it If the conditions of Theorem 2.1 are fulfilled on a countable set of intervals $[a_m;b_m], \ph m=1,2, \dots,$ and $\lim\limits_{m \to +\infty} a_m = + \infty$, then the system (1.1) is oscillatory. In this case the condition  $B(t) \ge 0$ may not be fulfilled outside of the set $\cup_{m=1}^{+\infty} [a_m; b_m]$.}

If the system (2.2) is oscillatory then  from the Sturm type comparison Theorem 3.8 of work [5] (see [5], p. 1511) it follows that for any $T \ge t_0$ there exists $T_1 > T$ such that the system (2.2) is oscillatory on $[T;T_1]$. Due to Remark 2.1 from here and from Theorem~ 2.1 we immediately get:

{\bf Corollary 2.1.} {\it Let the following conditions be satisfied:

\noindent
1') \ph $B(t) \ge 0, \ph t \ge t_0;$

\noindent
2') \ph there exists a continuously differentiable hermitian matrix - function $S(t)$  on $[t_0;+\infty)$ such that $A^*(t) + S(t) B(t) \in \Omega _n, \ph t\ge t_0$;

\noindent
3') the scalar system
$$
\sist{\phi' = \phh \sigma_S(t) \phi \ph +\ph  \frac{\lambda(B(t))}{n} \psi;}{\psi' = - tr(D_S(t)) \phi\ph - \ph \sigma_S(t) \psi, \ph t\ge t_0,} \eqno (2.3)
$$
is oscillatory.

\noindent
Then the system (1.1) is also oscillatory.} $\Box$

Let $\mu(t)$ be a real valued continuous function on $[t_0;+\infty)$. Consider the matrix equation
$$
B(t) X + X B(t) = 2\mu(t) I - A(t) - A^*(t), \phh t\ge t_0. \eqno (2.4)
$$
This equation has a solution on $[a;b]$, if in particular $B(t)$ and $-B(t)$ have no common eigenvalue for all $t\in [a;b]$ (e. g. $B(t) > 0, \ph t\in [a;b]$) (see [6], p. 207). One can easily show that if $B(t) \ge 0$ and $rank B(t) \ge n - 1$ for $t\in [a;b] \ph (t\ge t_0)$ then there exists some  real valued continuous function $\mu(t)$ on $[a;b] \ph ([t_0;+\infty))$ such that Eq. (2.4) has a solution on $[a;b]$ (on $[t_0;+\infty)$). It is not difficult to verify that if $X(t)$ is a solution of Eq. (2.4) then $S(t) \equiv \frac{X(t) + X^*(t)}{2}$ is a hermitian solution of    Eq. (2.4).

{\bf Theorem 2.2.} {\it Let the following conditions be satisfied:

\noindent
1) $B(t) \ge 0, \ph t\in [a;b]$;

\noindent
4) Eq. (2.4) has a continuously differentiable hermitian solution $S(t)$ on $[a;b]$;

\noindent
5) the scalar system
$$
\sist{\phi' = \phh \mu(t) \phi \ph +\ph  \frac{\lambda(B(t))}{n} \psi;}{\psi' = - tr(D_S(t)) \phi\ph - \ph \mu(t) \psi, \ph t\in [a;b],} \eqno (2.5)
$$
is oscillatory on $[a;b]$.

\noindent
Then the system (1.1) is also oscillatory on $[a,b]$.} $\Box$

Similar to Corollary 2.1 from here we obtain:

\pagebreak

{\bf Corollary 2.2.} {\it Let the following conditions be satisfied:

\noindent
1') $B(t) \ge 0, \ph t\ge t_0$;

\noindent
4') Eq. (2.4) has a continuously differentiable hermitian solution $S(t)$ on $[t_0;+\infty)$;

\noindent
5') the scalar system
$$
\sist{\phi' = \phh \mu(t) \phi \ph +\ph  \frac{\lambda(B(t))}{n} \psi;}{\psi' = - tr(D_S(t)) \phi\ph - \ph \mu(t) \psi, \ph t\ge t_0,} \eqno (2.6)
$$
is oscillatory.

\noindent
Then the system (1.1) is also oscillatory.} $\Box$

Obviously Theorem 2.1 and Theorem 2.2 as well as Corollary 2.1 and Corollary 2.2 are conditional results in that oscillation of the systems (2.2), (2.3), (2.5), (2.6) is only supposed rather than proved. The first of the following two assertions weakens the conditional character of Theorem 2.1 and Theorem 2.2 and the second one weakens the conditional character of Corollary 2.1 and Corollary 2.2. Set: $E(t) \equiv a_{11}(t) - a_{22}(t), \ph t\ge t_0$.

{\bf Theorem 2.3.} {\it Let the following conditions be satisfied:

\noindent
6) $a_{12}(t) \ge 0, \ph t\in [a;b];$

\noindent
7) $\int\limits_a^b \min\biggl[a_{12}(t)\exp\bigl\{-\int\limits_a^tE(\tau) d\tau\bigr\}, - a_{21}(t)\exp\bigl\{\int\limits_a^tE(\tau) d\tau\bigr\}\biggr]d t \ge \pi.$

\noindent
Then the system (2.1) is oscillatory on $[a;b]$.} $\Box$

{\bf Theorem 2.4.} {\it Let the following conditions be satisfied:

\noindent
6') $a_{12}(t) \ge 0, \ph t\ge t_0$;

\noindent
8) $\int\limits_{t_0}^{+\infty} a_{12}(t)\exp\bigl\{-\int\limits_{t_0}^tE(\tau) d\tau\bigr\}= - \int\limits_{t_0}^{+\infty}a_{21}(t)\exp\bigl\{\int\limits_{t_0}^tE(\tau) d\tau\bigr\} d t = +\infty$.

\noindent
Then the system (2.1) is oscillatory.} $\Box$

{\bf Remark 2.2}. {\it Theorem 2.4 is a generalization of the Leighton's oscillatory criterion (see [7], p. 70, Theorem 2.24).}

{\bf Remark 2.3.} {\it Another oscillatory criteria for the system (2.1) are proved in [5], which are applicable to the systems (2.3) and (2.6).}

Hereafter in this section we will assume that $B(t) \ge 0, \ph t\ge t_0,$ and $\sqrt{B(t)}$ is continu-\linebreak ously differentiable on $[t_0;+\infty)$. Consider the matrix equation
$$
\sqrt{B(t)} A^*(t) - \sqrt{B(t)}' = [\sqrt{B(t)} A^*(t) - \sqrt{B(t)}'] X \sqrt{B(t)}, \phh t\ge t_0.  \eqno (2.7)
$$
This equation has a solution on $[a;b]$ (on $[t_0;+\infty)$) if in particular $A^*(t) = \matr{A_1(t) \ph 0}{A_2(t) \ph 0}, \linebreak  B(t) = \matr{B_1(t) \ph 0}{\ph 0 \phh   0}, \ph t\ge t_0,$ where $A_1(t)$ and $B_1(t)$ are some matrices of dimension $m\times~ m \phantom{a}  (m < n)$ and $\det B_1(t) \ne 0, \ph t\in [a;b] \ph (t\ge t_0)$. In this case \linebreak $X(t) \equiv \matr{\sqrt{B_1(t)}^{-1} \ph 0}{ \phh 0 \phh  \phantom{aa} 0}, \ph t\in [a;b] \ph (t\ge t_0)$, is a solution of Eq (2.7) on $[a;b]$ (on $[t_0;+\infty)$).

Let $F(t)$ be a continuous matrix - function of dimension $n \times n$ on $[t_0;+\infty)$. Denote :

$$
\mathcal{D}_F(t) \equiv - \biggl[\frac{[\sqrt{B(t)} A^*(t) - \sqrt{B(t)}'] F(t) + F^*(t)[A(t)\sqrt{B(t)} - \sqrt{B(t)}']}{2}\biggr]' - \phantom{aaaaaaaaaaaaaaa}
$$
$$
\phantom{aaaa}- \biggl[\frac{[\sqrt{B(t)} A^*(t) - \sqrt{B(t)}'] F(t) + F^*(t)[A(t)\sqrt{B(t)} - \sqrt{B(t)}']}{2}\biggr]^2 - B(t) C(t), \ph t\ge t_0.
$$

{\bf Theorem 2.5}. {\it Let the following conditions be satisfied:

\noindent
1) $B(t) \ge 0, \ph t\in [a;b]$;

\noindent
9) Eq. (2.7) has a solution $F(t)$ on $[a;b]$ such that $[\sqrt{B(t)} A^*(t) - \sqrt{B(t)}'] F(t) +\linebreak \phh + F^*(t)[A(t)\sqrt{B(t)} - \sqrt{B(t)}']$ is continuously differentiable on $[a;b]$;

\noindent
10) the scalar equation
$$
\phi'' + \frac{tr(\mathcal{D}_F(t))}{n} \phi = 0, \phh t\in [a;b], \eqno (2.8)
$$
is oscillatory on $[a;b]$.

\noindent
Then the system (1.1) is also oscillatory on $[a;b]$.} $\Box$.

Taking into account Remark 2.1 on the strength of the Sturm's comparison theorem (see [8], p. 334, Theorem 3.1) from here we immediately get:

{\bf Corollary 2.3}. {\it Let the following conditions be satisfied:

\noindent
1') $B(t) \ge 0, \ph  t\ge t_0$;

\noindent
9') Eq. (2.7) has a solution $F(t)$ on $[t_0;+\infty)$ such that $[\sqrt{B(t)} A^*(t) - \sqrt{B(t)}'] F(t)$ is continuously differentiable on   $[t_0;+\infty)$;

\noindent
10') the scalar equation
$$
\phi'' + \frac{tr(\mathcal{D}_F(t))}{n} \phi = 0, \phh t\ge t_0,
$$
is oscillatory.

\noindent
Then the system (1.1) is also oscillatory.} $\Box$.

{\bf Remark 2.4.} {\it If $B(t) > 0, \ph t\in [a;b] \ph (t\ge t_0)$ and $\sqrt{B(t)} A^*(t) - \sqrt{B(t)}'\not\equiv 0$ then $F(t) \equiv \sqrt{B^{-1}(t)}, \ph t\in [a;b] \ph (t\ge t_0)$ is the unique solution of Eq. (2.7), and if in addition $\sqrt{B(t)}$ and $A(t)$ are permutable (e. g. $A(t) = \sum\limits_{j=1}^N \alpha_j(t) K^j(t), \ph \sqrt{B(t)}~ =\linebreak =\sum\limits_{j=1}^N \beta_j(t) K^j(t), \ph t\ge t_0$, where $\alpha_j(t) ,\ph  \beta_j(t), \ph j=\overline{1,N},$ are some continuous functions on $[t_0;+\infty)$, $K(t)$ is a continuous square matrix function on $[t_0;+\infty)$; more detailed information about permutable matrices one can find in [6 pp. 199 - 207]) then it can be shown that
$$
tr(\mathcal{D}_{\sqrt{B^{-1}}}(t)) = - tr \biggl[\biggl(\frac{A(t) + A^*(t)}{2}\biggr)' + \biggl(\frac{A(t) + A^*(t)}{2}\biggr)^2 + B(t) C(t) - (\sqrt{B(t)}'\sqrt{B^{-1}(t)})' +
$$
$$
\phantom{aaaaaaaaaaaaaaaaaaaaaaaaaaaa} + \frac{1}{2}\Bigl(A(t) + A^*(t)\Bigr)\sqrt{B(t)}'\sqrt{B^{-1}(t)}\biggr], \ph t\in [a;b],\ph (t\ge t_0).
$$}

{\bf 2.2. Examples}. In this section we present some examples  demonstrating the capacities of the obtained results.

Example 2.1. Consider the matrix equation
$$
\Phi'' + K(t)\Phi = 0, \phh t\ge t_0 > 0, \eqno (2.10)
$$
$$
\mb{where} K(t) \equiv \matrr{a_1\sin \mu_1 t  +a_2 \sin \mu_2 t \phh  \frac{b\cos \mu_3 t}{t^\alpha} \phh  \phh \phh \phh 0}{\phh \frac{b\cos \mu_3 t}{t^\alpha} \phh a_1\sin \mu_1 t  +a_2 \sin \mu_2 t  \phh \phh \frac{c\sin \mu_4 t}{t^\beta}}{\phh \ph 0 \phh \phh \phh \phh \frac{c\sin \mu_4 t}{t^\alpha} \phh \phh a_1\sin \mu_1 t  +a_2 \sin \mu_2 t }, \phantom{aaaaaaaaaaaaaaaa}
$$
$a_1, \ph a_2, \ph \alpha, \ph \beta, \ph b, \ph c, \ph \mu_j, \ph j= \overline{1,4},$ are some real constants, $a_j \ne 0, \ph j=1,2, \ph \alpha >~ 1, \linebreak \beta >~ 1, \ph \mu_1 / \mu_2$ is irrational. This Equation is equivalent to the system (1.1) for \linebreak $A(t) \equiv 0, \ph B(t) \equiv I$, where $I$ is the identity matrix of dimension $3\times 3, \ph C(t) \equiv - K(t)$. Therefore according to Theorem 2.1 Eq. (2.10) is oscillatory provided the scalar system
$$
\sist{\phi' = \phantom{aaaaaaaaaaaaaaaaaaaaaaaaaaa} \frac{1}{3} \psi;}{\psi' = - 3(a_1\sin\mu_1 t + a_2 \sin \mu_2(t)) \phi, \phantom{aaaaaaaa} t\ge t_0,}
$$
is oscillatory, which is equivalent to the oscillation  of the scalar equation
$$
\phi'' + (a_1\sin \mu_1 t + a_2 \sin \mu_2 t) \phi = 0, \phh t\ge t_0.
$$
This equation is oscillatory (see [9], Corollary 1). Therefore the last system also is \linebreak oscillatory. From here it follows that Eq. (2.10) is oscillatory. The eigenvalues of the matrix $K(t)$ are equal $\lambda_\pm(t) \equiv a_1\sin \mu_1 t + a_2 \sin \mu_2 t \pm \sqrt{\frac{b^2\cos^2 \mu_3 t}{t^{2\alpha}} + \frac{c^2\sin^2 \mu_4 t}{t^{2\beta}}}, \ph \lambda_1(t) \equiv a_1\sin \mu_1 t + a_2 \sin \mu_2 t, \ph t\ge t_0.$ This shows that the Theorems 5, 6 of work [10], the Theorems 1, 2, 3 of work [11] are not applicable to Eq. (2.10). The conditions of the remaining results of these works and the conditions of the results of the works  [1 - 4, 12 -~ 14] contain arbitrary parameter-functions. Therefore it is very difficult  to guess the applicability of these results to Eq. (2.10).

Example 2.2.  Set: $A^*_1(t) \equiv \matr{\cos t \phh a(t)}{-\overline{a(t)} \ph \cos t}, \phh \phh B_1(t) \equiv \frac{1}{t} \matr{1 \ph \sin t}{\sin t \ph 1}, \linebreak C_1(t) \equiv \matr{-1/t + \alpha \cos t \phh c(t)}{ \phh \overline{c(t)} \phh \phantom{aaaaa} \beta \sin t}, \ph t\ge 1,$ where $\alpha, \beta \in (-\infty; +\infty), \ph a(t)$  and $c(t)$ are some continuous functions on $[1;+\infty)$. Consider the system
$$
\sist{\Phi' = A_1(t)\Phi + B_1(t)  \Psi;}{\Psi' = C_1(t) \Phi - A_1^*(t) \Psi, \ph t\ge 1,} \eqno (2.11)
$$
We will use Corollary 2.1 to show that this system is oscillatory. Set:   $S(t) \equiv 0, \ph t\ge 1$. Then it is not difficult to verify that $A^*_1(t)+ S(t)B(t) \in \Omega_n, \ph  \lambda(B_1(t))= \frac{1 - |\sin t|}{t}, \ph t\ge~ 1$. $\sigma_S(t) = \cos t, \ph D_S(t) = - C_1(t), \ph t\ge~ 1$ and $\int\limits_1^{+\infty} \frac{1 - |\sin t|}{2 t} \exp\biggl\{ -2 \int\limits_1^t \cos \tau d\tau\biggr\} d t =\linebreak  = \int\limits_1^{+\infty}\biggl( \frac{1}{t}  - \alpha \cos t - \beta \sin t\biggr)\exp\biggl\{2 \int\limits_1^t \cos \tau d\tau\biggr\} d t =  + \infty$. By Theorem 2.4 from here it follows that all conditions of Corollary 2.1 for the system (2.11) are fulfilled. Therefore the system (2.11) is oscillatory.

Example 2.3. Consider the system
$$
\sist{\Phi' = \phantom{aaaa} K_2(t) \Psi;}{\Psi' = - K_2(t) \Phi, \phantom{aaaa} t \ge 0,} \eqno (2.12)
$$
where $K_2(t) \equiv diag\{\nu \sin t, \dots , \nu \sin t\}, \ph t \ge 0, \ph \nu \ge \frac{\pi}{2}$. Obviously Corollary 2.1 is not applicable to this system (the condition 1) is not fulfilled).
Note that By Theorem 2.3  for this system for $a= 2\pi m, \ph b = \pi (2 m +1)$ the conditions of Theorem 2.1 are satisfied for each $m=1,2, \dots.$ Due to Remark 2.1 from here it follows that the system (2.12) is oscillatory.

{\bf Remark 2.5.} {\it No result of works  [1- 4, 12 - 14]    is applicable to the systems \linebreak (2.10) - (2.12)}.

{\bf Remark 2.6}. {\it Suppose $A(t) \equiv 0, \ph B(t) = - C(t) \equiv I, \ph t \ge 0$, where $I$ is the identity matrix. It is evident that in this case for the system (1.1) the conditions 1) - 3) of Theorem 2.1 are fulfilled on the arbitrary interval $[a;b] (\subset [0;+\infty))$ and the condition 4) is fulfilled only if $b - a \ge \pi$. It also is evident that for this case $(\Phi_0(t), \Psi_0(t))$, where \linebreak $\Phi_0(t) \equiv diag\{\sin t, \dots, \sin t\}, \ph \Psi_0(t) \equiv diag \{\cos t, \dots, \cos t\}$, is a prepared solution to the system (1.1). This solution is not oscillatory on $[\varepsilon; \pi - \varepsilon]$ for each $\varepsilon \in (0;1)$.  Therefore in the inequality 7) we may not replace $\pi$ by a number less than $\pi$ (in this sense the condition 7) is sharp).}

Example 2.4. Set $\mathcal{M}(t) \equiv \max\{ \sin t, 0\}, \ph A_2(t) \equiv \matrr{\mu(t) \ph 2\sin t \phantom{aa} (1 + \mathcal{M}(t)) \cos t} {\ph  0 \phh \mu(t) \phh (1 + \mathcal{M}(t)) \sin t}{ \ph 0 \phh 0 \phh \phh \phh \mu(t)}, \linebreak B_2(t) \equiv \matrr{1 \phh 0 \phh 0}{0 \phh 1 \phh 0}{ 0 \phh 0 \phh \mathcal{M}(t)}, \ph C_2(t) \equiv \matrr{-\mathcal{M}(t) \sin^2 t \phh 0 \phh \ph  0}{\phh \ph 0  \phantom{aaaaaa} -1 \phh \ph 0}{\phh \ph 0 \phantom{aaaaaaaa} 0 \ph -2 \mathcal{M}(t)}, \ph t\ge t_0,$ where $\mu(t)$ is a continuous real valued function on $[t_0;+\infty)$. Consider the system
$$
\sist{\Phi' = A_2(t)\Phi + B_2(t)  \Psi;}{\Psi' = C_2(t) \Phi - A_2^*(t) \Psi, \ph t\ge 1,} \eqno (2.13)
$$
 One can readily check that  $\ph S(t) \equiv \matrr{0 \ph - \sin t \ph - \cos t}{- \sin t \ph 0 \ph - \sin t}{- \cos t - \cos t \ph  0}$ is a solution to the matrix equation
$$
B_2(t) X + X B_2(t) = 2\mu(t) I - A_2(t) - A_2^*(t), \ph t\ge t_0.
$$
After some simple calculations we get: $\lambda(B_2(t)) = \mathcal{M}(t), \ph tr D_S(t) = - 2 \sin^2 t, \ph t\ge~ t_0$.
On the basis of Corollary 2.2 and Theorem 2.4  from here we  conclude that if the function $\mathcal{G}(t) \equiv \il{t_0}{t}\mu(\tau) d\tau, \ph t\ge t_0,$ is bounded then the system (2.13) is oscillatory.

Example 2.5. Let $a_{jk}(t), \ph c_{jk}(t), \ph j,k = \overline{1,3}, \ph  \beta(t)$ be continuous functions on $[t_0;+\infty)$ such that $a_{12}(t) + a_{22}(t) = a_{11}(t) + a_{21}(t), \ph a_{31}(t) = a_{32}(t), \ph c_{jk}(t) = \overline{c_{kj}(t)}, \ph j,k = \overline{1,3}, \ph \beta(t) \ge 0, \ph t\ge t_0; \ph \beta(t)$ and $ Re [a_{11}(t) + a_{21}(t) + a_{33}(t)] - \frac{\beta'(t)}{2\beta(t)}$ is continuously differentiable on $[t_0;+\infty)$. Set; $A^*_3(t) \equiv (a_{jk}(t))_{j,k = 1}^3, \ph  C_3(t)\equiv (c_{jk}(t))_{j,k = 1}^3, \linebreak B_3(t) \equiv \matrr{1 \ph 1 \phh 0}{ 1 \ph 1\phh 0}{ 0\ph 0 \ph \beta(t)}, \ph t\ge t_0$. Consider the system
$$
\sist{\Phi' = A_3(t) \Phi + B_3(t) \Psi;}{\Psi' = C_3(t) \Phi - A_3^*(t) \Psi, \phh t\ge t_0.} \eqno (2.14)
$$
It is not difficult to verify that $B_3(t) \ge 0, \ph \sqrt{B_3(t)} = \matrr{\frac{\sqrt{2}}{2} \ph \frac{\sqrt{2}}{2} \phh 0}{\frac{\sqrt{2}}{2} \ph \frac{\sqrt{2}}{2} \phh 0}{\ph 0\phantom{aa} 0 \phantom{aa} \sqrt{\beta(t)}}, \ph t\ge t_0$, and the matrix function  $F_3(t) \equiv \matrr{\frac{\sqrt{2}}{4}  \ph \frac{\sqrt{2}}{4}  \phantom{aa} 0}{\frac{\sqrt{2}}{4}  \ph \frac{\sqrt{2}}{4}  \phantom{aa} 0}{\ph 0 \phantom{aa} 0 \ph \frac{1}{\sqrt{\beta(t)}}}, \ph t\ge t_0,$ is a solution to the matrix equation
$$
\sqrt{B_3(t)} A^*_3(t) - \sqrt{B_3(t)}' = [\sqrt{B_3(t)} A_3(t) - \sqrt{B_3(t)}'] X \sqrt{B_3(t)}
$$
on $[t_0;+\infty)$. So we have all data for calculation of  $tr (\mathcal{D}_{F_3}(t)), \ph t\ge t_0$. After some simple arithmetic operations we obtain
$$
tr (\mathcal{D}_{F_3}(t)) = - \biggl[ Re (a_{11}(t) + a_{21}(t) + a_{33}(t)) - \frac{\beta'(t)}{2\beta(t)}\biggr]' - \biggl[ Re (a_{11}(t) + a_{21}(t) + a_{33}(t)) - \frac{\beta'(t)}{2\beta(t)}\biggr]^2  -
$$
$$
\phantom{aaaaaaaaaaaaaaaaaaaaaaaaaaaaaaaaaaaaaaaaaaaa} - c_{11}(t)  - 2 Re c_{12}(t) - c_{22}(t), \phh t\ge t_0.
$$
By Theorem 2.5 (Corollary 2.3) the system (2.14) is oscillatory on $[a;b]$ (is oscillatory) provided the scalar equation $$
\phi'' + \frac{tr(\mathcal{D}_{F_3}(t))}{3} \phi = 0, \phh t\ge t_0,
$$
is oscillatory on $[a;b]$ (is oscillatory).

{\bf Remark 2.7}. {\it One can readily check that $\lambda(B_3(t)) \equiv 0, \ph t\ge t_0$. Therefore Theorem~ 2.1 and Theorem 2.2 as well as Corollary 2.1 and Corollary 2.2 are not applicable to the system~ (2.14).}

\vskip 20pt

\centerline{\bf \S 3. Proof of the main results}.

 \vskip 20pt

{\bf 3.1.  Auxiliary propositions}.
 Let $f(t), \ph g(t), \ph h(t) \ph f_1(t), \ph g_1(t), \ph h_1(t)$  be real valued continuous functions on $[t_0;+\infty)$. Consider the Riccati equations:
$$
y' + f(t) y^2 + g(t) y + h(t) = 0, \phh t \in [a;b]; \eqno (3.1)
$$
$$
y' + f_1(t) y^2 + g_1(t) y + h_1(t) = 0, \phh t \in [a;b]; \eqno (3.2)
$$
and the following inequalities
$$
\eta' + f(t) \eta^2 + g(t) \eta + h(t) \ge 0, \phh t \in [a;b]; \eqno (3.3)
$$
$$
\eta' + f_1(t) \eta^2 + g_1(t) \eta + h_1(t) \ge 0, \phh t \in [a;b]; \eqno (3.4)
$$

{\bf Remark 3.1}. {\it If $f(t) \ge 0$ \ph for $t\in [a;b]$ then every solution of the linear equation
$$
y' + g(t) y + h(t) = 0
$$
on the interval $[a;b]$ is also a solution of the inequality (3.3).}

{\bf Remark 3.2.} {\it Every solution of Eq. (3.2) on the interval $[a;b]$ is also a solution of the inequality (3.4).}

{\bf Theorem 3.1}. {\it Let Eq. (3.2) has a real solution $y_1(t)$ on $[a;b]$, and let the following conditions be satisfied: $f(t) \ge 0$ and $\int \limits_a^t\exp\biggl\{\int\limits_a^\tau [f(s)(\eta_0(s) + \eta_1(s)) + g(s)] d s\biggr\} \biggl[(f_1(\tau) - f(\tau)) y_1^2(\tau) + (g_1(\tau) - g_1(\tau) - g(\tau)) y_1(\tau) + h_1(\tau) - h(\tau)\biggr] d \tau \ge 0, \ph t\in [a;b],$ where $\eta_0(t)$ and $\eta_1(t)$ are solutions of the inequalities (3.3) and (3.4) respectively on $[a;b]$ such that $\eta_j(a) \ge y_1(a), \ph j=1,2.$. Then for every $\gamma_0 \ge Y_1(a)$ Eq. (3.1) has a real valued solution $y_0(t)$ on $[a;b]$, satisfying the initial condition $y_0(a) = \gamma_0$.}

\noindent
Proof. By analogy of the proof of Theorem 3.1 from [15].

Consider the inequality
$$
y' + f(t) y^2 + g(t) y + h(t) \le  0, \phh t \in [a;b]; \eqno (3.5)
$$

{\bf Lemma 3.1.} {\it If $f(t) \ge 0, \ph t\in [a;b]$ then Eq. (3.1) has a solution on [a;b] if and only if the inequality (3.5) has a solution on $[a;b]$.}

Proof. Obviously every solution of Eq. (3.1) is also a solution  of the inequality (3.5). Let $y_1(t)$ be a solution to the inequality (3.5) on $[a;b]$. Set $\widetilde{h}(t) \equiv - (y_1'(t) + f(t) y_1^2(t) + g(t) y_1(t)), \ph t\in [a;b]$. Since
$
y_1(t)' + f(t) y_1^2(t) + g(t) y_1(t) + \widetilde{h}(t) \le  0, \phh t \in [a;b],
$
we have
$$
h(t) \le \widetilde{h} (t), \phh t\in [a;b]. \eqno (3.6)
$$
Consider the equation
$$
y' + f(t) y^2 + g(t) y + \widetilde{h}(t) =  0, \phh t \in [a;b].
$$
Using Theorem 3.1 to this equation and Eq. (3.1) and taking into account (3.6) we conclude that Eq. (3.1) has a solution on [a;b]. The lemma is proved.

{\bf Lemma 3.2.} {\it For any two square matrices $M_1\equiv (m_{ij}^1), \ph M_2\equiv (m_{ij}^2)$ the equality
$$
tr (M_1 M_2) = tr (M_2 M_1)
$$
is valid.}

Proof. We have $tr (M_1 M_2) = \sum\limits_{j=1}^n(\sum\limits_{k=1}^n m_{jk}^1 m_{kj}^2) = \sum\limits_{k=1}^n(\sum\limits_{j=1}^n m_{jk}^1 m_{kj}^2) = \sum\limits_{k=1}^n(\sum\limits_{j=1}^n m_{kj}^2 m_{jk}^1) = tr (M_2 M_1).$ The lemma is proved.

{\bf Lemma 3.3.} {\it Let the following conditions be satisfied:

\noindent
1$^*$) $f(t) \ge 0, \ph t\ge t_0$; \ph 2$^*$) $h(t) \ge 0, \ph t\ge t_0$; \ph 3$^*$) $\int\limits_{t_0}^{+\infty} f(\tau) \exp\biggl\{- \int\limits_{t_0}^\tau g(s) d s\biggr\} d\tau = + \infty$; \linebreak

\noindent
4$^*$) For some $t_1 \ge t_0$ the equation
$$
y' + f(t) y^2 + g(t) y + h(t) = 0, \phh t\ge t_0  \eqno (3.7)
$$
has a real valued solution on $[t_1;+\infty)$.

\noindent
Then Eq. (3.7) has a positive solution on $[t_1;+\infty)$.}

Proof. Let (according to the condition 4$^*$))  the function  $y(t)$ be a  solution of Eq. (3.7) on $[t_1;+\infty)$  for some $t_1 \ge t_0$ and let $y_1(t)$ be another solution of Eq. (3.7) with
$$
y_1(t_1) > y(t_1). \eqno (3.8)
$$
Then from 1$^*$) it follows that $y_1(t)$ exists on $[t_1;+\infty)$ (see [16]). Show that
$$
y_1(t) > 0, \phh t\ge t_1. \eqno (3.9)
$$
Suppose for some $t_2 \ge t_1$
$$
y_1(_2) \le 0.  \eqno (3.10)
$$
Consider the linear equation
$$
y' + \xi(t) y + h(t) = 0, \phh t\ge t_2,
$$
where $\xi(t) \equiv f(t) y_1(t)  + h(t), \ph t\ge t_2$. Obviously $y_1(t)$  is a solution of this equation. Then by Cauchy's  formula we have:
$$
y_1(t) = \exp\biggl\{-\int\limits_{t_2}^t \xi(\tau) d \tau\biggr\}\biggl[y_1(t_2) - \int\limits_{t_2}^t\exp\biggl\{\int\limits_{t_2}^\tau \xi(s) d s\biggr\} h(\tau) d \tau\biggr], \phh t\ge t_2.
$$
From here from 2$^*$) and  (3.10) it follows that
$$
y_1(t) \le 0, \phh t\ge t_2. \eqno (3.11)
$$
From 3$^*$) and from the easily verifiable equality
$$
\int\limits_{t_2}^{+\infty} f(\tau) \exp\biggl\{- \int\limits_{t_2}^\tau g(s) d s\biggr\} d\tau = \exp\biggl\{- \int\limits_{t_0}^{t_2} g(s) d s\biggr\}\times \phantom{aaaaaaaaaaaaaaaaaaaaaaaaaaaaaaaaaaaaaa}
$$
$$
\phantom{aaaaaaaaaaaaaaaaaaaaa} \times \biggl[\int\limits_{t_0}^{+\infty} f(\tau) \exp\biggl\{- \int\limits_{t_0}^\tau g(s) d s\biggr\} d\tau  - \int\limits_{t_0}^{t_2} f(\tau) \exp\biggl\{- \int\limits_{t_0}^\tau g(s) d s\biggr\} d\tau  \biggr]
$$
it follows that $\int\limits_{t_2}^{+\infty} f(\tau) \exp\biggl\{- \int\limits_{t_2}^\tau g(s) d s\biggr\} d\tau = +\infty$.
From here from 1$^*$) and (3.11) it follows
$$
\nu_{y_1}(t_2) \equiv \int \limits_{t_2}^{+\infty} f(\tau) \exp\biggl\{ - \int\limits_{t_2}^\tau[2 f(s) y_1(s) = g(s)] d s\biggr\} = +\infty.
$$
On the other hand from1$^*$) and (3.8) it follows that (see [16]). $\nu_{y_1}(t_2) < + \infty$. The obtained contradiction proves (3.9). The lemma is proved.

{\bf 3.1.  Proof of Theorem 2.1}. Let $(\Phi(t), \Psi(t))$ be a prepared solution of the system (1.1). Show that $det \Phi(t)$ vanishes on $[a;b]$. Suppose that it is not true. Then  \hskip2pt $det \Phi(t) \ne~0, \linebreak t \in [a;b]$.  It follows from here that the hermitian matrix $Y_0(t) \equiv \Psi(t) \Phi ^{-1}(t), \ph t\in [a;b],$ is a solution to the matrix Riccati equation
$$
Y' + Y B(t) Y + A^*(t) Y + Y A(t) - C(t) = 0, \phh t\in[a;b]. \eqno (3.12)
$$
Let $S(t)$ satisfies the condition 2). In (3.12) make the substitution: $Y = Z + S(t), \ph t\ge T$. We obtain
$$
Z' + Z B(t) Z + [A^*(t) + S(t) B(t)] Z + Z [A(t) + B(t) S(t)] + D_S(t) = 0, \phh t \in [a;b].
$$
Obviously the hermitian matrix - function $Z_0(t) \equiv Y_0(t) + S(t), \ph t \in [a;b],$ is a solution to this equation on $[a;b]$. Therefore
$$
[tr Z_0(t)]' + \frac{\lambda(B(t))}{n} [tr Z_0(t)]^2 +\phantom{aaaaaaaaaaaaaaaaaaaaaaaaaaaaaaaaaaaaaaaaaaaaaaaaaaaa}
$$
$$
\phantom{aaaa} + tr [(A^*(t) + S(t) B(t)) Z_0(t) + Z_0(t) (A(t) + B(t) S(t))] +  tr D_S(t) \le 0, \phh t \ge T. \eqno (3.13)
$$
Since $A(t) + S(t) B(t) \in \Omega_n, \ph t\ge t_0$, it is not difficult to verify that $tr [(A(t) + \linebreak S(t) B(t)) Z_0(t) + Z_0(t) (A^*(t) + B(t) S(t))] = 2\sigma_S(t) tr Z_0(t), \ph t \ge T.$ From here and from (3.13) we get:
$$
[tr Z_0(t)]' + \frac{\lambda(B(t))}{n} [tr Z_0(t)]^2 + 2 \sigma_S(t) Z_0(t) + tr D_S(t) \le 0, \phh t\ge T. \eqno (3.14)
$$
By 1) we have $\lambda(B(t)) \ge 0, \ph t\in [a;b]$. By virtue of Lemma 3.1. from here and from (3.14) it follows that the equation
$$
y' + \frac{\lambda(B(t))}{n} y^2 + 2\sigma_S(t) y  + tr D_S(t) = 0, \phh t\in [a;b],
$$
has a solution $y(t)$ on $[a;b]$.
Therefore (see [17]) the functions
$$
\phi(t) \equiv \exp\biggl\{\il{a}{t}\biggl[\frac{\lambda(B(\tau))}{n} y(\tau) + \sigma_S(\tau)\biggr] d \tau\biggr\}, \phh \psi(t) \equiv y(t) \phi(t), \phh t \in [a;b],
$$
form a non oscillatory solution $(\phi(t), \psi(t))$ of the system (2.2) on $[a;b]$. Hence the system (2.2) is not oscillatory on $[a;b]$, which contradicts the condition 3) of the theorem. The obtained contradiction completes the proof of the theorem.

{\bf 3.2. Proof of Theorem 2.2}. Suppose the system (1.1) is not oscillatory on $[a;b]$. Then there exists a prepared solution $(\Phi(t), \Psi(t))$ of the system (1.1) such that $det \Phi(t) \ne~ 0, \linebreak t\in [a;b]$. Then for the hermitian matrix - function $Y(t) \equiv \Psi(t) \Phi^{-1}(t), \ph t\in [a;b]$, the following equality takes place
$$
Y'(t) + Y(t) B(t) Y(t) + A^*(t) Y(t) + Y(t) A(t) - C(t) = 0, \phh t\in [a;b].
$$
In this equality make the substitution $Y(t) \equiv U(t) + S(t), \ph t\in [a;b]$, where $S(t)$ is a hermitian solution of Eq. (2.3) on $[a;b]$. Taking into account 4) we get:
$$
U'(t) + U(t) B(t) U(t) + (2 \mu(t)I - L(t)) U(t) + U(t) L(t) + D_S(t) = 0, \ph t\in [a;b], \eqno (3.15)
$$
where $L(t)\equiv B(t) S(t) + A(t), \ph t\in [a;b]$. By Lemma 3.2 \ph $tr[L(t)U(t) - U(t)L(t)] \equiv ~0, \linebreak t\in [a;b]$. From here and from (3.15) we obtain
$$
[tr(U(t))' + \frac{\lambda(B(t))}{n}[tr (U(t))]^2 + 2\mu(t) tr (U(t))  + tr(D_S(t)) \le 0, \ph t\in [a;b].
$$
By virtue of Lemma 3.1 from here it follows that the equation
$$
y' + \frac{\lambda(B(t))}{n} y^2 + 2\mu (t) y + tr(D_S(t)) = 0, \phh t\in [a;b],
$$
has a solution $y_1(t)$ on [a;b].
 Therefore  the functions
$$
\phi_1(t) \equiv \exp\biggl\{\il{a}{t}\biggl[\frac{\lambda(B(\tau))}{n} y_1(\tau) + \mu(\tau)\biggr] d \tau\biggr\}, \phh \psi_1(t) \equiv y(t) \phi_1(t), \phh t \in [a;b],
$$
form a non oscillatory solution $(\phi_1(t), \psi_1(t))$ of the system (2.4) on $[a;b]$. Therefore the system (2.4) is not oscillatory on $[a;b]$, which contradicts the condition 4) of the theorem.  The obtained contradiction proves the theorem.

{\bf3.3.  Proof of Theorem 2.3}.  In the system (2.1) make the substitutions:
$$
\sist{\phi = \exp\biggl\{\il{a}{t}a_{11}(\tau) (\tau) d \tau\biggr\} \rho \sin \theta;}{\psi = \exp\biggl\{ \il{a}{t}a_{22}(\tau) d \tau\biggr\} \rho \cos \theta, \phh t \in [a;b].} \eqno (3.16)
$$
We will get:
$$
\sist{\rho' \sin \theta + \theta' \rho  \cos \theta = A_{12}(t) \rho \cos \theta;}{\rho' \cos \theta - \theta' \rho  \sin \theta = A_{21}(t) \rho \sin \theta, \phh t\in [a;b],} \eqno (3.17)
$$
where $A_{12}(t) \equiv a_{12}\exp\biggl\{- \il{a}{t}E(\tau) d \tau\biggr\}, \ph A_{21}(t) \equiv a_{21}(t)\exp\biggl\{ \il{a}{t}E(\tau) d \tau\biggr\}, \ph t\in [a;b].$ This system is equivalent to the system (2.1) in the sense that to each nontrivial solution $(\phi(t), \psi(t))$ of the system (2.1) corresponds the solution $(\rho(t), \theta(t))$ of the system (3.17) with $\rho(t) > 0, \ph t\in [a;b)$, defined by (3.16). Let us multiply the first equation of the system (3.17) by $\cos \theta$ and the second one by $\sin \theta$ and subtract from the first obtained the second one. We  get:
$$
\theta' \rho = \rho [A_{12}(t) \cos^2 \theta - A_{21}(t) \sin ^2 \theta], \phh t\in [a;b]. \eqno (3.18)
$$
Let $(\phi_0(t), \psi_0(t))$ be a nontrivial solution of the system (2.1) and let $(\rho_0(t), \theta_0(t))$ be the solution of the system (3.17) corresponding to $(\phi_0(t), \psi_0(t))$. Then $\rho_0(t) \ne 0, \ph t\in [a;b]$, and therefore by (3.18) the following equality takes place
$$
\theta_0'(t) = A_{12}(t) \cos^2 \theta_0(t) - A_{21}(t) \sin^2 \theta_0(t) = \frac{1}{2}\Bigl[A_{12}(t) - A_{21}(t) + (A_{12}(t) + A_{21}(t)) \cos 2\theta_0(t)\Bigr],
$$
$t\in [a;b].$ From here it follows
$$
\theta_0'(t) \ge  \frac{1}{2} \Bigl[A_{12}(t) - A_{21}(t) - |A_{12}(t) + A_{21}(t)|\Bigr] =  \min\{A_{12}(t), - A_{21}(t)\}, \ph t\in [a;b].
$$
Let us integrate this inequality from $a$ to $b$. Taking into account the conditions of the theorem we  get:
$$
\theta_0(b) - \theta_0(a) \ge  \il{a}{b}\min\Bigl\{A_{12}(\tau), - A_{21}(\tau)\Bigr\} d\tau \ge \pi.
$$
Due to (3.16) from here it follows that $\phi_0(t)$ has at least one zero on $[a;b]$. The theorem is proved.

{\bf 3.4. Proof of Theorem 2.4}. Suppose the system (2.1) is not oscillatory. Then (see [17]) the equation                                                $$
y' + a_{12}(t) y^2 + E(t) y - a_{21}(t)  = 0, \phh t\ge t_0, \eqno (3.19)
$$
has a solution $y(t)$ on $[t_1;+\infty)$ for some $t_1 \ge t_0$.
Set: $u(t) \equiv a_{12}(t) \exp\biggl\{- \il{t_1}{t} E(\tau) d\tau\biggr\}, \linebreak  w(t) \equiv - a_{21}(t) \exp\biggl\{\int\limits_{t_1}^t E(\tau) d \tau\biggr\}, \ph t\ge t_1.$ In Eq. (3.19) make the substitution
$$
y = z \exp\biggl\{-  \il{t_1}{t}E(\tau) d \tau\biggr\}, \phh t \ge t_1.
$$
We obtain
$$
z' + u(t) z^2 + w(t) = 0, \phh t\ge t_1, \eqno (3.20)
$$
 Show that
$$
\ilp{t_1} u(\tau)\exp\biggl\{\il{t_1}{t}2 u(\tau) d\tau\il{t_1}{\tau} w(s) sd s\biggr\} d t = + \infty. \eqno (3.21)
$$
By 8) we have $\il{t_1}{t} w(\tau) d \tau = - \il{t_1}{t} a_{21}(\tau)\exp\biggl\{- \il{t_1}{\tau} E(s) d s\biggr\} d \tau \ge 0, \ph t \ge t_2,$ for some $t_2 \ge t_1$.
From here and from 8) it follows (3.21).
In Eq. (3.20) make the substitution \linebreak $z = U - \il{t_1}{t} w(\tau) d \tau, \ph t \ge t_1$. We  get:
$$
U' + u(t) U^2 - 2u(t) \il{t_1}{t} w(\tau) d \tau  \hskip 2pt U + u(t) \biggl[\il{t_1}{t} w(\tau) d \tau\biggr]^2  = 0, \phh t\ge t_1. \eqno (3.22)
$$
Since   Eq. (3.19) has a real valued solution on $[t_1;+\infty)$, from the
 substitutions of dependent variables, made  above,  it can be seen that Eq. (3.22) has a real valued solution on $[t_1;+\infty)$. On the strength of  Lemma 3.3 from here from (3.21) and from the inequalities $a_{12}(t) \ge~ 0, \linebreak u(t)\biggl[\il{t_1}{t} w(\tau) d\tau\biggr]^2 \ge~ 0, \ph t \ge t_1,$ it follows that Eq. (3.22) ha a positive solution  $U_0(t)$ on $[t_1;+\infty)$. Then $Z_0(t) \equiv U_0(t) - \il{t_1}{t} w(\tau) d\tau, \ph t\ge t_1,$ is a solution to Eq. (3.20) on $[t_1;+\infty)$ such that
$$
Z_0(t) > - \il{t_1}{t} w(\tau) d\tau, \phh t\ge t_1. \eqno (3.23)
$$
It follows from (3.20) that
$$
Z_0(t) = Z_0(t_1) - \il{t_1}{t} u(\tau) Z_0^2(\tau)d\tau - \il{t_1}{t} w(\tau) d \tau, \phh t\ge t_1. \eqno (3.24)
$$
From here and from (3.23) it follows that
$$
0 \le \il{t_1}{t} u(\tau) Z_0^2(\tau)d\tau < Z_0(t_1), \phh t\ge t_1 \eqno (3.25)
$$
$(Z_0(t_1) = U_0(t_1) > 0)$. Taking into account 8) from here we get: \linebreak $\biggl[Z_0(t_1) - \il{t_1}{t} u(\tau) Z_0^2(\tau)d\tau - \il{t_1}{t} w(\tau)\biggr]^2 \ge 1, \ph t \ge T,$ for some $T\ge t_1$. From here and from (3.24) it follows that $Z_0^2(t) \ge 1, \ph t \ge T.$ Therefore by 8) we have $\ilp{T}u(\tau) Z_0^2(\tau) d\tau \ge \ilp{T}u(\tau) d\tau = + \infty,$ which contradicts (3.25). The obtained contradiction completes the proof of the theorem.

{\bf 3.5. Proof of Theorem 2.5}. Suppose for some prepared solution $(\Phi(t), \Psi(t))$ of the system (1.1) $\det \Phi(t) \ne 0, \ph t\in[a;b].$ Then for the hermitian matrix - function $Y(t) \equiv~ \Psi(t) \Phi^{-1}(t), \ph t\in [a;b]$, the following equality holds
$$
Y'(t) + Y(t) B(t) Y(t) + A^*(t) Y(t) + Y(t) A(t) - C(t) = 0, \phh t\in [a;b].
$$
Multiplying  both sides of this equality at left and at right by $\sqrt{B(t)}$, and Taking into account the  equality $(\sqrt{B(t)} Y(t) \sqrt{B(t)})' = \sqrt{B(t)} Y'(t) \sqrt{B(t)} +  \sqrt{B(t)}' Y(t) \sqrt{B(t)} +\linebreak + \sqrt{B(t)} Y(t) \sqrt{B(t)}', \ph t\in [a;b],$ we get:
$$
Z'(t) + Z^2(t) + [\sqrt{B(t)} A(t) - \sqrt{B(t)}'] Y(t) \sqrt{B(t)} + \phantom{aaaaaaaaaaaaaaaaaaaaaaaaaaaaa}
$$
$$
\phantom{aaaaaaaaaaaaa}+ \sqrt{B(t)} Y(t) [A^*(t) \sqrt{B(t)} - \sqrt{B(t)}'] - \sqrt{B(t)} C(t) \sqrt{B(t)} = 0, \phh t\in [a;b],
$$
where $Z(t) \equiv \sqrt{B(t)} Y(t) \sqrt{B(t)}, \ph t\in [a;b]$. From here and from the condition 2) we obtain
$$
Z'(t) + Z^2(t) + L(t) Z(t) + Z(t) L^*(t) - \sqrt{B(t)} C(t) \sqrt{B(t)} = 0, \phh t\in [a;b].
$$
where $L(t) \equiv [\sqrt{B(t)} A^*(t) - \sqrt{B(t)}'] F(t), \ph t\in [a;b].$ Making substitution \linebreak $Z(t) \equiv V(t) - \frac{L(t) + L^*(t)}{2}, \ph t\in [a;b],$ in this equality we get:
$$
V'(t) + V^2(t) + \mathcal{D}_F(t) + \frac{L(t) - L^*(t)}{2} V(t) - V(t) \frac{L(t) - L^*(t)}{2} +  \phantom{aaaaaaaaaaaaaaaaaaa}
$$
$$
\phantom{aaaaaaaaaaaaaaaaaaaaaaaaaa}+ B(t) C(t) - \sqrt{B(t)} C(t) \sqrt{B(t)} = 0, \phh t\in [a;b]. \eqno (3.26)
$$
By Lemma 3.2 $tr\biggl[\frac{L(t) - L^*(t)}{2} V(t) - V(t) \frac{L(t) - L^*(t)}{2}\biggr]  = tr [B(t) C(t) - \sqrt{B(t)} C(t) \sqrt{B(t)}] \equiv~ 0, \linebreak t\in [a;b].$ From here and from (3.26) we obtain:
$$
[tr (V(t)]' + \frac{1}{n}[tr (V(t))]^2 + tr (\mathcal{D}_F(t)) \le 0, \phh t\in [a;b].
$$
By Lemma 3.1 from here it follows that the Riccati equation
$$
y' + \frac{1}{n} y^2 + tr (\mathcal{D}_F(t)) =0, \phh t\in [a;b],
$$
has a solution $y(t)$ on $[a;b]$. Then the function $\phi(t) \equiv \exp\biggl\{\il{a}{t} \frac{y(\tau)}{n} d \tau\biggr\}, \ph t\in [a;b],$ is a non vanishing  solution of Eq. (2.8)  on $[a;b]$. Therefore Eq. (2.8) is not oscillatory on $[a;b]$, which contradicts the condition   10) of the theorem. The obtained contradiction completes the proof of the theorem.

\vskip 20 pt

\centerline{\bf References}

\vskip 20 pt

\noindent
1. L. Li, F. Meng and Z. Zheng, Oscillation Results Related to Integral Averaging Technique\linebreak \phantom{a} for Linear Hamiltonian Systems, Dynamic Systems and Applications 18 (2009), \ph \linebreak \phantom{a} pp. 725 - 736.

\noindent
2. F. Meng and  A. B. Mingarelli, Oscillation of Linear Hamiltonian Systems, Proc. Amer.\linebreak \phantom{a} Math. Soc. Vol. 131, Num. 3, 2002, pp. 897 - 904.

\noindent
3. Q. Yang, R. Mathsen and S. Zhu, Oscillation Theorems for Self-Adjoint Matrix \linebreak \phantom{a}   Hamiltonian
 Systems. J. Diff. Equ., 19 (2003), pp. 306 - 329.

\noindent
4. Z. Zheng and S. Zhu. Hartman Type Oscillatory Criteria for Linear Matrix Hamiltonian  \linebreak \phantom{a} Systems. Dynamic  Systems and Applications, 17 (2008), pp. 85 - 96.

\noindent
5.  G. A. Grigorian. Oscillatory Criteria for the Systems of Two First - Order Linear\linebreak \phantom{a} Ordinary Differential Equations. Rocky Mountain Journal of Mathematics, vol. 47,\linebreak \phantom{a} Num. 5, 2017, pp. 1497 - 1524

\noindent
6.  F. R. Gantmacher, Theory of Matrix. Second Edition (in Russian). Moskow,,\linebreak \phantom{a} ''Nauka'', 1966.

\noindent
7. C. A. Swanson, Comparison and Oscillation Theory of Linear Differential Equations.\linebreak \phantom{a} Academic Press. New York and London, 1968.

\noindent
8. Ph. Hartman, Ordinary differential equations, Second edition, The Jhon Hopkins \linebreak  \phantom{aa}University, Baltimore, Maryland, 1982.

\noindent
9. G. A. Grigorian, On one Oscillatory Criterion for The Second Order Linear
 Ordinary \linebreak \phantom{a} Differential Equations. Opuscula Math. 36, Num. 5 (2016), 589–601. \\ \phantom{a}
   http://dx.doi.org/10.7494/OpMath.2016.36.5.589

\noindent
10. L. H. Erbe, Q. Kong and Sh. Ruan, Kamenev Type Theorems for Second Order Matrix\linebreak \phantom{aa}  Differential Systems. Proc. Amer. Math. Soc. Vol. 117, Num. 4, 1993, 957 - 962.

\noindent
11. R. Byers, B. J. Harris and M. K. Kwong, Weighted Means and Oscillation Conditions\linebreak \phantom{a}  for Second Order Matrix Differential Equations. Journal of Differential Equations\linebreak \phantom{a} 61, 164 - 177 (1986).

\noindent
12. G. J. Butler, L. H. Erbe and A. B. Mingarelli, Riccati Techniques and Variational\linebreak \phantom{aa} Principles in Oscillation Theory for Linear Systems, Trans. Amer. Math. Soc. Vol. 303,\linebreak \phantom{aa} Num. 1, 1987, 263 - 282.

\noindent
13. A. B. Mingarelli, On a Conjecture for Oscillation of Second Order Ordinary Differential\linebreak \phantom{aa} Systems, Proc. Amer. Math. Soc., Vol. 82. Num. 4, 1981, 593 - 598.

\noindent
14. Q. Wang, Oscillation Criteria for Second Order Matrix Differential Systems Proc.\linebreak \phantom{aa} Amer.  Math. Soc. Vol. 131, Num. 3, 2002, 897 - 904.

\noindent
15. G. A. Grigorian.  On Two Comparison Tests for Second-Order Linear  Ordinary \linebreak \phantom{aa}   Differential Equations (Russian) Differ. Uravn. 47 (2011), no. 9, 1225 - 1240; translation \linebreak \phantom{aa}in Differ.
 Equ. 47 (2011), no. 9 1237 - 1252, 34C10.

\noindent
16. G. A. Grigorian, Properties of Solutions of Riccati Equation. Journal of Contemporary \linebreak  \phantom{aa}Mathematical Analysis. 2007, vol. 42, No. 4, pp.184 - 197.

\noindent
17. G. A. Grigorian, On the Stability of Systems of Two First - Order Linear Ordinary\linebreak \phantom{aa} Differential Equations, Differ. Uravn., 2015, vol. 51, no. 3, pp. 283 - 292.

\end{document}